\documentclass[12pt,reqno]{amsart}
\usepackage{amsfonts,latexsym,amssymb,amsmath,amsthm}
\usepackage{graphicx,hyperref,cite,subfigure}

\begin{document}

\title{Architectural form as space-time cell}
\author{Luisa Consiglieri, Vitor Consiglieri}
\address{Luisa Consiglieri, Independent Researcher Professor, European Union}
\urladdr{\href{http://sites.google.com/site/luisaconsiglieri}{http://sites.google.com/site/luisaconsiglieri}}
\address{Vitor Consiglieri,
Architect and Emeritus Professor, Faculty of Architecture, Technical University of Lisbon, Portugal}

\begin{abstract}
The architecture has its basis in a dialectic search of new choices of representation. 
 We deal with the form on the contemporary architecture 
under two approaches: expression and content.
 We examine how mathematical principles based on natural growth can be applied in architectural design in order to create a
 dynamic, rather than static, structure. The dynamic process of a cell 
{\it and its growth} provides the basic structure. 
We exemplify the impact of the new forms on the new society that already began.
\end{abstract}

\keywords{mapping (transformation); architectural form;
shape; continuity; growth}
\subjclass[2010]{00A67; 97G40; 97I60; 97M80}

\maketitle

\section{Introduction}
In a contemporary analysis, the gestalts of the beginning of the XX century follow a society interested in the environmental endeavors and find the respect with the climatic system and the natural resources of the earth. Moreover, in the period of the German Bauhaus School at the beginning of XX century, Paul Klee defended that the study of the form should reflect the theory of Everything\footnote{The unification of quantum mechanics and relativity as an extension of the classical mechanics.}
 and it should introduce the space-time concept. Some innovative architectural designs and theories are seeking out a form best suited to conveying the spatiotemporal experiences phrased by Relativity 
\cite{hat}. Different forms can be obtained by the mathematical mappings
\cite{cc06,cc09},
 whence the importance of the mathematics research in architecture is highlighted.

Indeed the dialectics between nature and society has being analyzed
 by several authors 
\cite{ank,bech,eld,hva,mur}.
{\it If we reduce the art to the exact and automatic repetition of the Nature and confine it only to revive and remember an already felt emotion and a subconscious thing, the signification of the art would be mediocre or null 
\cite[p. 126]{sal}. An architectural concept depends upon two factors: the existence of a certain number of correlated elements and of one condition of relationship 
\cite[p. 90]{sal}. It creates in the human being a new emotion without which would not be art and would confine to repeat and to revive an already felt emotion 
\cite[p. 126]{sal}.}

The conceptual art appears in the architecture as a healthy interaction between the urban life and nature. The phenomenon of the random urbanism is surpassed, and the eco-cities are powered by renewable energy and free from pollution, and provide an opportunity for sustainable urban design 
\cite{cas,wat,yea}.  We deal with the form
 in relation to society and the issue of sustainable development, respectively,
\begin{itemize}
  \item   the society of the individual as a single or familiar element;
  \item  and the correction of climatic environment of the earth.
\end{itemize}

The conscious processes of composing architectural objects by organic forms should take into account how powerful is the graphical representation of the mathematical mappings.
Our  main concern is that
the principal question is not just about making beautiful buildings but also about how the architectural forms should be the key element of a sustainable development work on the urban context.
To this end, we begin by adapting architectural form to an environmental context.
We focus on the fundamental role played by
the expression and the content on architectural forms.

\section{The expression of the form}
\label{ef}

The crucial quality of shape, no matter of what kind, lies in its organization, and when we think of it in this way we call it form 
\cite{alex}. The expression of a form results from the search on the possibility in realizing ample spaces with less supports. The aesthetic value of the object's composition  relies on
 the continuity of the so called tension lines that shape the fusion of interior and exterior spaces. Spaces can be designed to encourage interaction or seclusion (see \cite{alh} and references therein).
 The relationship between the interior and the exterior, or the psychological weight that each one of these spaces has through the living conditions and the human activities, is given by the equilibrium of the structural components. Before reaching the curvilinear topological stage, the conception of the structure was a passive system corresponding to the walls that discharge their load in the soil. For instance, when the distance between the supports is small, the radius of curvature is found such that the forces are compensated. When the distance between the supports is increased, the stresses affect the vault, changing its configuration. Then, the expression is based on the intercommunication factors corresponding to the morphocontinuity
\cite{cc06,cc09}, in particular to the following meanings of continuity:
\begin{itemize}
  \item 	the interior/exterior continuity;
\item the light/shade continuity.
\end{itemize}

\begin{figure}
\begin{center}
\resizebox{1\hsize}{!}{\includegraphics{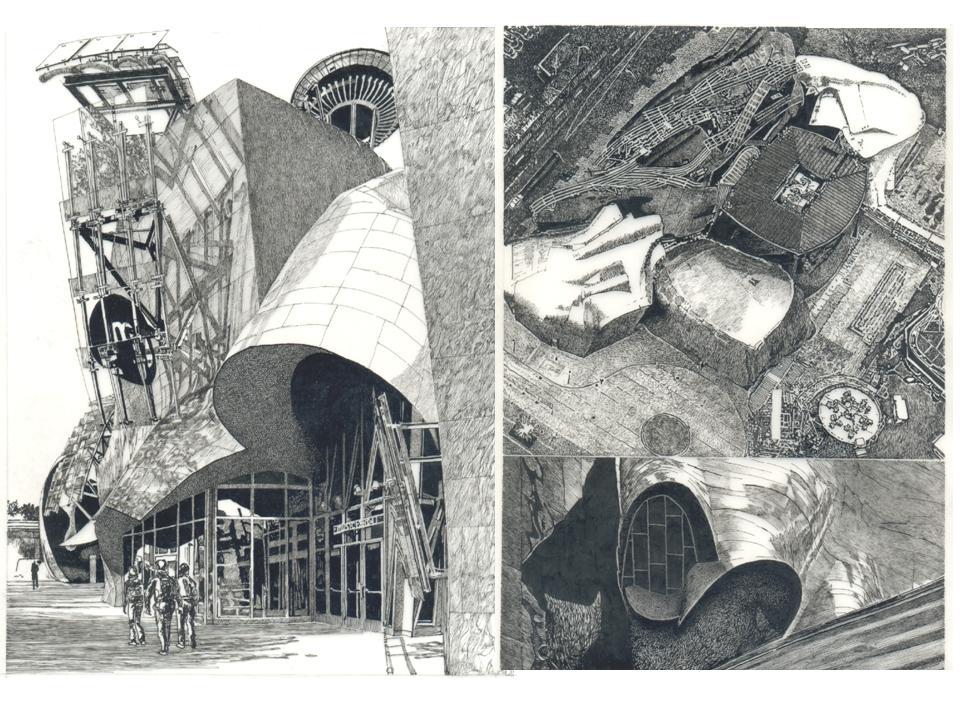}}
\caption{Experience Music Project (Seattle, 1995-2000) by Frank Gehry.}
\label{fig1}
\end{center}
\end{figure}

\begin{figure}
\begin{center}
\subfigure[Guggenheim Museum (Bilbao, 1997)]{\resizebox*{6cm}{!}{\includegraphics{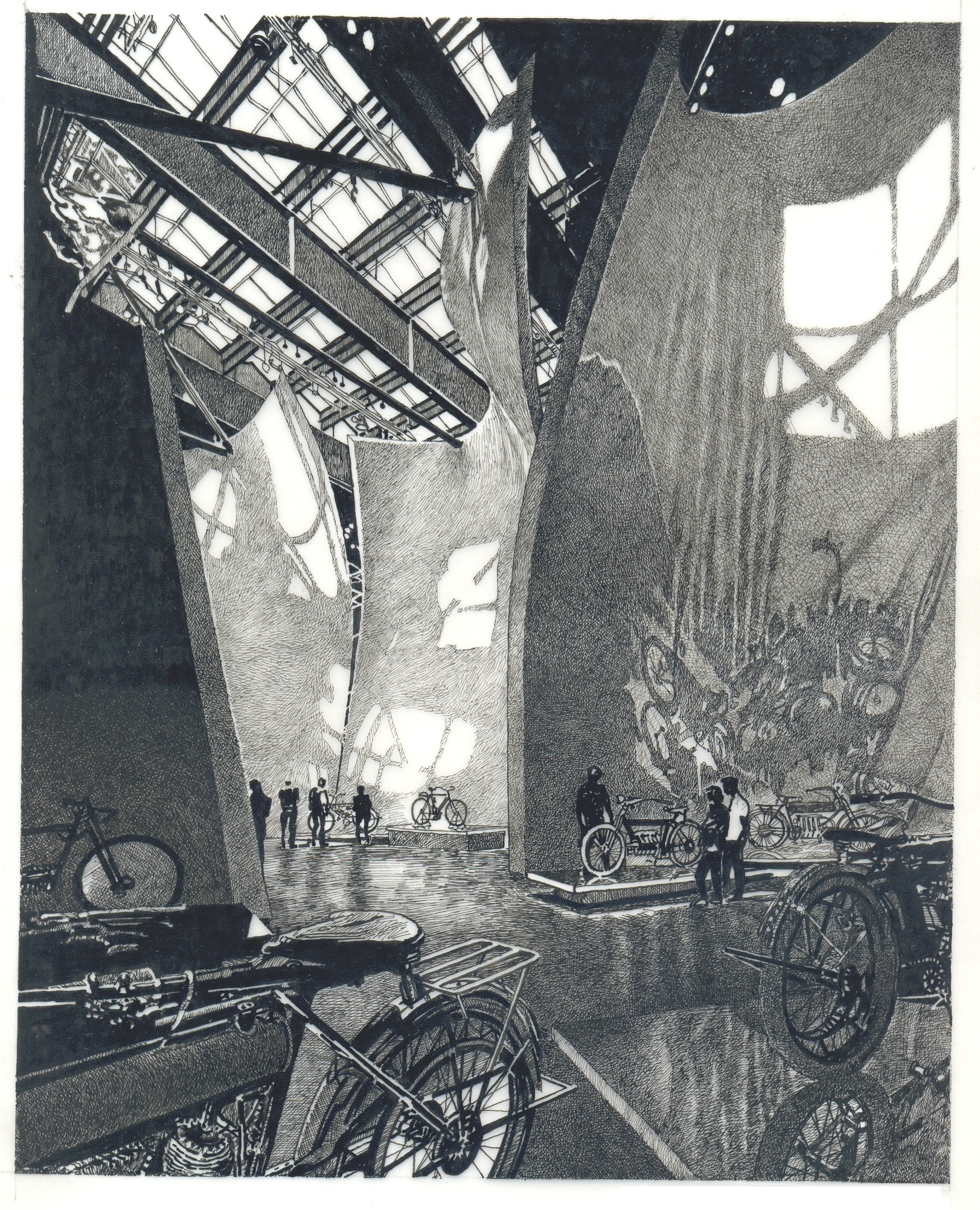}}}
\subfigure[Bank (Berlin, 1998-2000)]{\resizebox*{6cm}{!}{\includegraphics{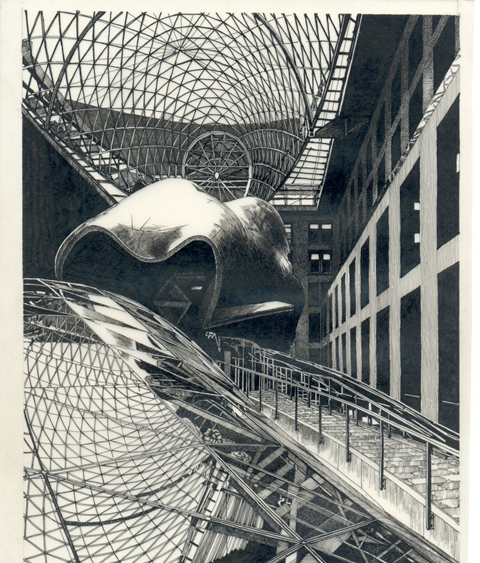}}}
\caption{Architect: Frank Gehry.}
\label{fig2}
\end{center}
\end{figure}

\begin{figure}
\begin{center}
\subfigure[Milwaukee Art Museum (Wisconsin, 2001) by Santiago Calatrava]{\resizebox*{6cm}{!}{\includegraphics{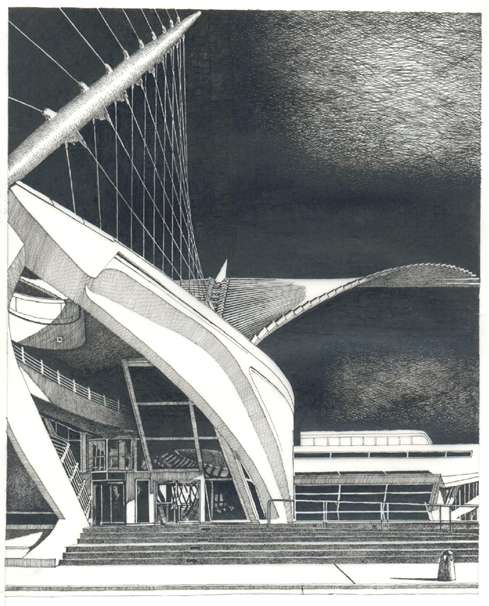}}}
\subfigure[Guggenheim Museum (Las Vegas) by Rem Koolhaas]{\resizebox*{6cm}{!}{\includegraphics{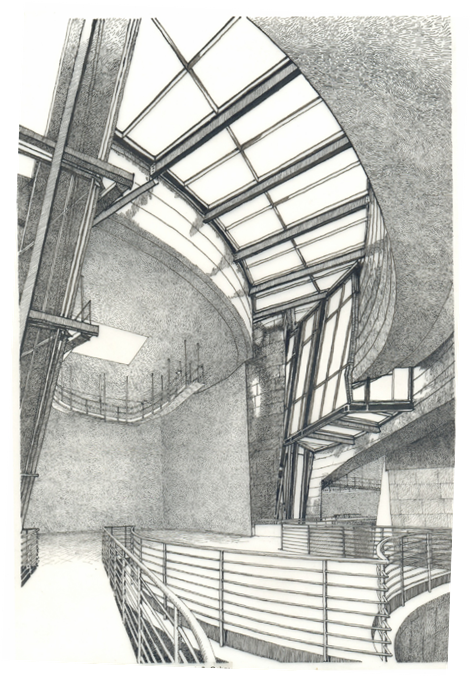}}}
\caption{}\label{fig3}
\end{center}
\end{figure}

In these interpretations, the connection between interior and exterior or between light and shade can be supplied by the glassy structures that give visual communicability or visibility. In particular, the environment in an inner space changes according to the transparencies given by the existence of the windows. The more or less size of the windows infers a graduation of light/shade effects. Also the existence of doors or arches contributes to that connection. Nowadays a greater connection between the inner and outer spaces comes by means of curvilinear forms, attracting the two spaces in permanent fusion (Figure \ref{fig1}).
 The passage from the classical systems into new forms, still represented by continuous mappings, follows such that it gives a visual dynamism and it supplies energetic exchange of thermodynamics, ventilation and fluid dynamics, and light and shading (Figure \ref{fig2}).
 Moreover, the topological configuration that integrates the building creates the connection between light and shade (Figure \ref{fig3}).
 The different situations of opening of the inner space compel considerable light/shade continuity.

The network of multifunctional open green spaces and human-made infrastructures can be finite or infinite. We understand these two classifications of the carry-over of the activity from the inner space into the outer space as: finite, when the outer space is a terrace with its bounds, that is, it is a space area of shared action; and infinite, when there does not exist boundaries in the outer space, that is, the space looses sight in the horizon given a panorama space or of contemplation. Indeed, it represents 
{\it a world inside the interior of the world}
 as argued by Louis Kahn \cite[p. 391]{nor}.

\section{The content of the form}
\label{cf}

The climatic changes are already deeply influencing the whole planet. This transformation obliges the search of processes that improve the 
ecological performance and, in particular, reduce the waste. 
In the architectonic processes, two main 
environmental factors are the production and the acclimation. 

\begin{figure}
\begin{center}
\resizebox{1\hsize}{!}{\includegraphics{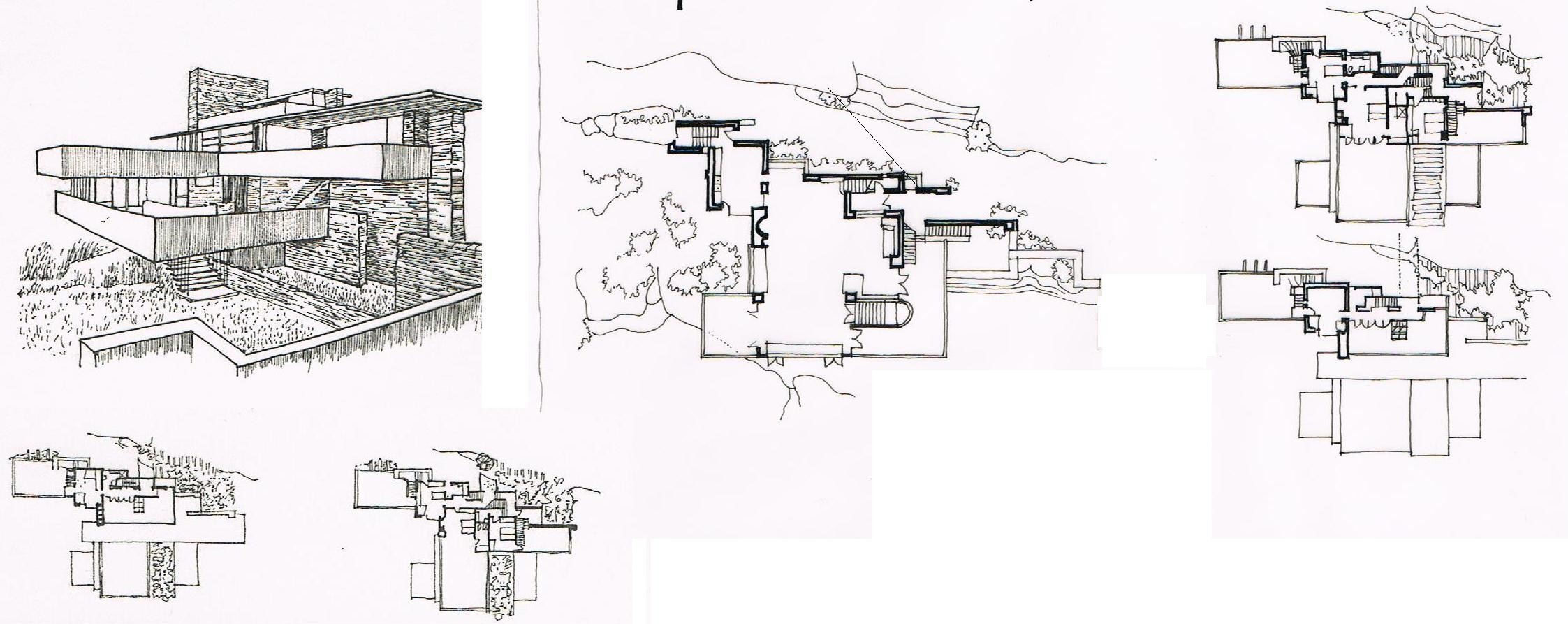}}
\caption{Fallingwater House in Pennsylvania by Frank L. Wright.}
\label{fig9}
\end{center}
\end{figure}

We begin by discuss how the production infer in the ecological performance. In the beginning of XX century, constructive stability and structural density were the main objectives.  For instance, the expression of the Fallingwater House 
(Figure \ref{fig9}), by Frank
L. Wright, is given by the general orthogonality of the structures coincident with the stability principles of the concrete-steel material. With the concrete shells, this orthogonal behavior disappears and a dynamical stage with effective 3D 
forms appears in constant change. Moreover, the Euclidean forms are substituted by non-Euclidean ones. These forms can be studied under the properties of mathematical mappings. Indeed, they are membranes of cells as defined in the morphologic rhythm 
\cite{cc06,cc09}. The form expresses harmony between the inner and the outer spaces as for instance the works of Saarinen 
\cite{cc10}.

\begin{figure}
\begin{center}
\resizebox{1\hsize}{!}{\includegraphics{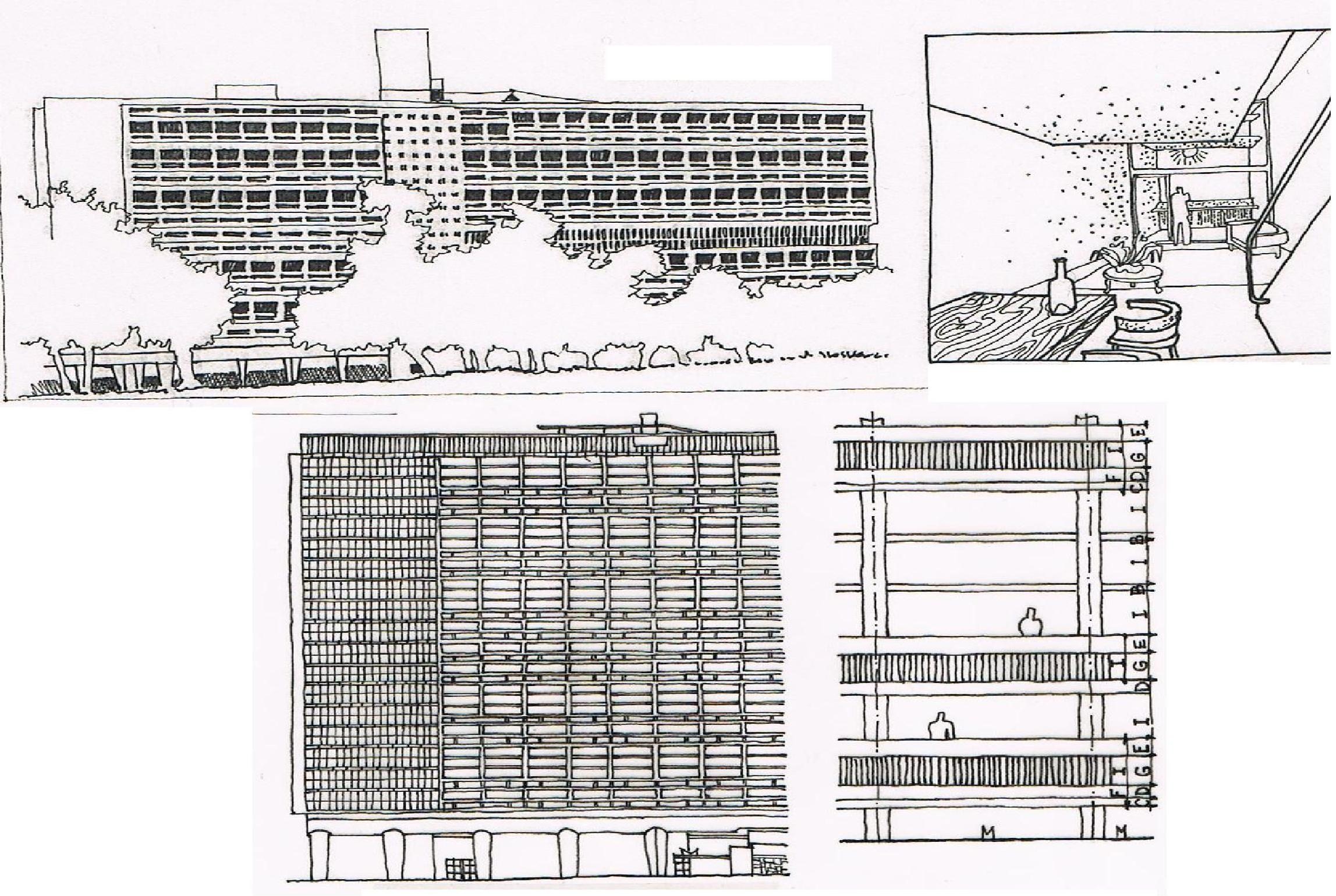}}
\caption{Unit\'e d'Habitation in Marseille, by Le Corbusier.}
\label{fig10}
\end{center}
\end{figure}

The second cause in the ecological performance is the divergence about the acclimation of the buildings, either housing or public, that is, the existence of natural and artificial acclimation. In the Le Corbusier's buildings, for instance the Unit\'e d'Habitation in Marseille 
(Figure \ref{fig10}), we find a concern in this sense. 
The study of light-breaker caused by vertical and horizontal surface elements creates a primordial design to the east, south and west facades. By that time, Mies van der Rohe stands up for the acclimation given by the electrical artificial system. He thought it as the more suitable reason for the internal environment and for inward way of living. Indeed, the artificial acclimated buildings lead the architects in concentrating only on the aesthetics of the form. The architects and engineers, such as Felix Candela, Pier L. Nervi, 
Eduardo Torroja, and Frei Otto, among others, stand up that the building structure must have a relation with the internal, visual, constructive and plastic form. The structure becomes the form itself. Buckminster Fuller laid the problem in the following manner: the society develops from the visible into the invisible and the architecture has to keep up with it. In fact, the electronic network transforms our world into a cybernetic one.

\begin{figure}
\begin{center}
\resizebox*{8cm}{!}{\includegraphics{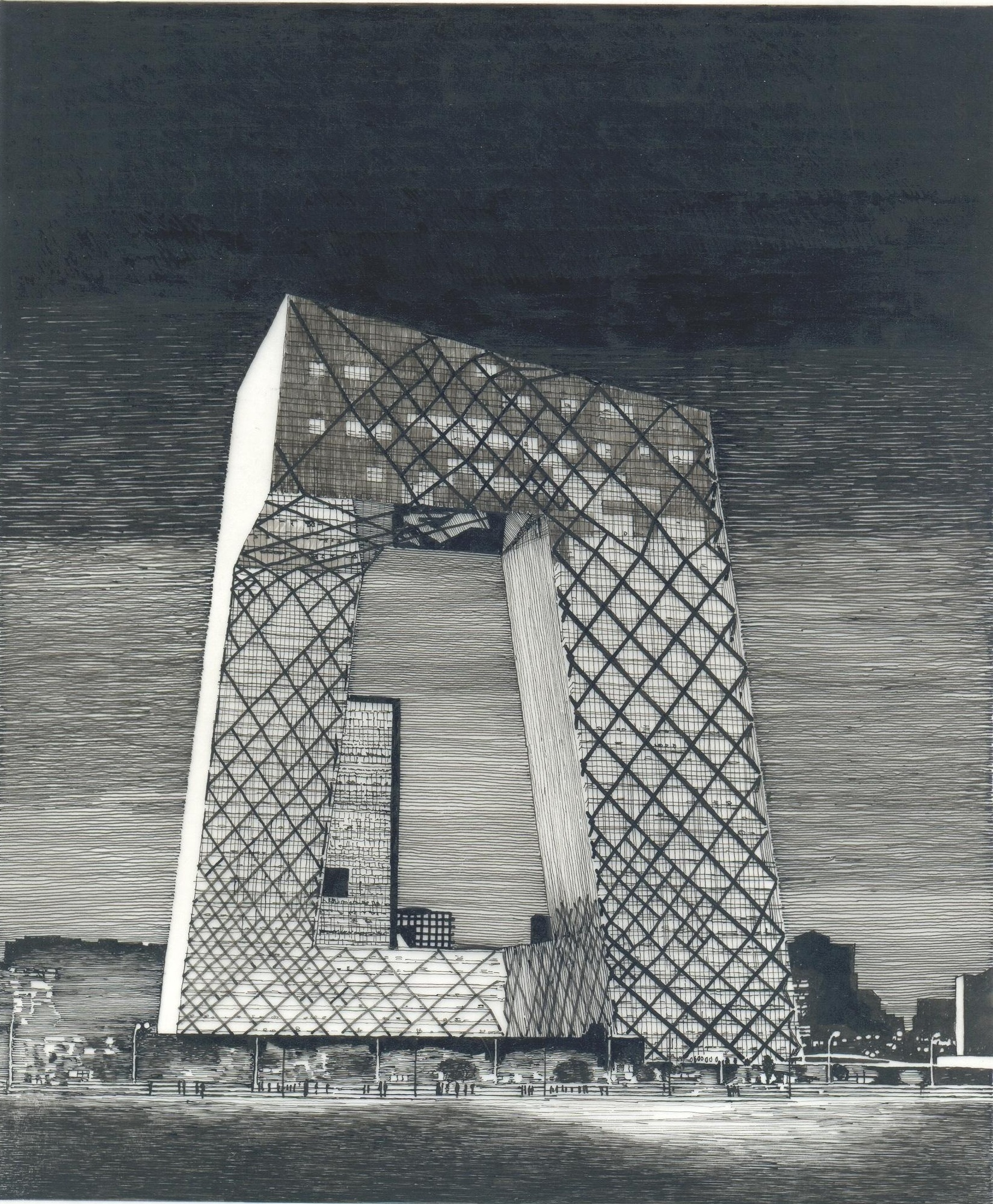}}
\caption{Central Chinese Television Tower (Beijing, 2004-08)
 by Rem Koolhaas and Ole Scheeren (architects) and Ove Arup 
\& Partners (engineers).}
\label{fig11}
\end{center}
\end{figure}

With the rectification of the ecological errors done to the global system, the study of new forms is being connected with new materials. The glassy structures by one hand valorize
 the volumetry and on the other hand encourage the inward well-being and comfort. Thus the urban environment is enriched by the glassy forms. Recently, the greatest attempt in serving the mechanical systems such as heating, ventilation and cooling is the looping configuration 
\cite{Ril}
(Figure \ref{fig11}). However, it has a great objection: big expenses. On the search on exuberance, it costs high energetic charges both in production and management.

At the present century, the interdisciplinary teams accomplished by the ecological concern deal the beauty of the form with the flora. The phenomena under the biological point of view are reduced to the hereditary properties of the organisms, and consequently to the study of populations. Mathematically, this study is commonly known as dynamical systems. In order to compare the biological phenomena with the architecture the study must not be reduced to the imitation of the nature but must eliminate the energetic excesses. The buildings' projects, beyond the structural systems that support the object, ought to include ventilation and shelter systems for climatic conditions. For instance, the Olympic Stadium of Beijing in China is built like a bird nest that lets in light and air, and at the same time it filtrates the cold and the wind. The structure protects both the athletes and the audience. Moreover, it brings together a good thermo installment in letting to pass in the solar light and a natural ventilation without any artificial methods and equipments. This complex system can be seen as a continuous form keeping its nature inspiration.
 The binary pairs fade on the process of continuity \cite{so,elo}.
 
\section{The morphocontinuity}
 
With the study of the efficacy of cooling/ventilation and insulation in conjunction with the structural stability of the buildings, the materials quality is explored in a biological way of looking. We enter in a new stage of the architectural relationship between the building and human living, notwithstanding Leonardo DaVinci had already said: In Nature, everything has always one reason. If we understand such reason, the experience is not anymore needed.
 The morphocontinuity stands for the intrinsic property of form as
a spatial cell that continuously depends on instants of time $t\in\mathbb{R}$
(for details, see \cite{cc06,cc09}). Here, we define {\it cell} by the set
\[
C=\{(x,y,z)\in D_f:\ f(x,y,z;t)\leq 1\},
\]
where the spatial and temporal variables ($x,y,z$ and $t$)
are assumed to be dimensionless, the origin stands for a
reference space-time point, and $f$ is a real function, continuous on the
time real line, but not necessarily
continuous on its spatial domain $D_f\subset \mathbb{R}^3$.

\begin{figure}
\begin{center}
%\begin{minipage}{120mm}
\subfigure[$t=1$]{\resizebox*{4cm}{!}{\includegraphics{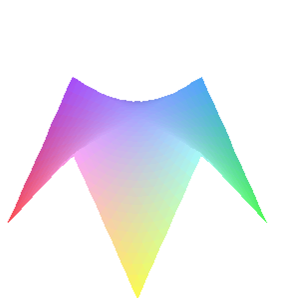}}}
\subfigure[$t=2$]{\resizebox*{4cm}{!}{\includegraphics{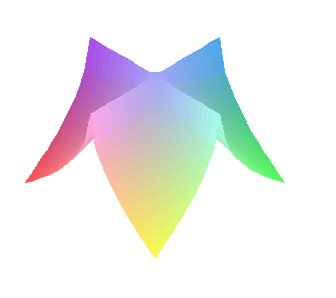}}}
\subfigure[$t=4$]{\resizebox*{4cm}{!}{\includegraphics{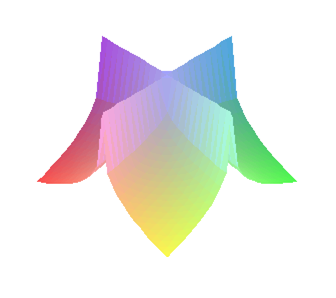}}}
\caption{The graphical representation of a form defined by the set
$\{(x,y,z)\in [-2,2]^2\times\mathbb{R}:$ $z=|xy|^{1/t}\}$,
 for the time parameter $t>0$.}
\label{fig4}
%\end{minipage}
\end{center}
\end{figure}

 When the research in providing the architectural objects against the climatic changes is care of, a form can be obtained through the mathematical mappings. Figure \ref{fig4}
 shows the spatiotemporal characterization of a form reflecting the opening behavior, for instance, to the Philips Pavilion (Brussels) by Le Corbusier and Xenakis. 
 The geometric volumes that till now had mostly served as motifs as architectural autonomous forms with a stability character, they are now cells composed of a most dynamic expression.

\begin{figure}
\begin{center}
\subfigure[Reichstag's Dome (Berlin, 1999)]{\includegraphics[width=5cm]{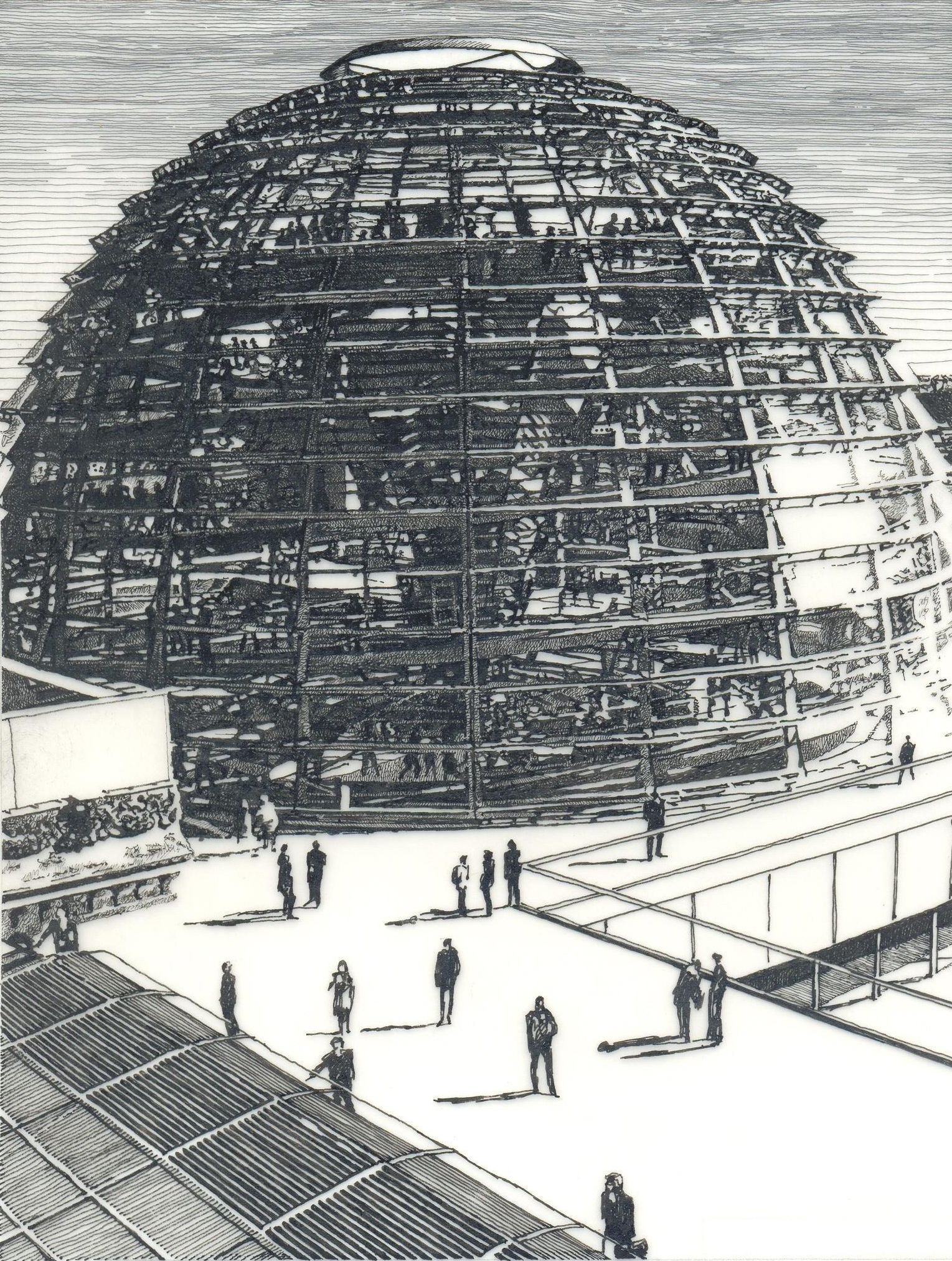}}
\subfigure[The Swiss Reinsurance Headquarters 30th St. Mary Axe (London, 1997-2004)]{\includegraphics[width=5cm]{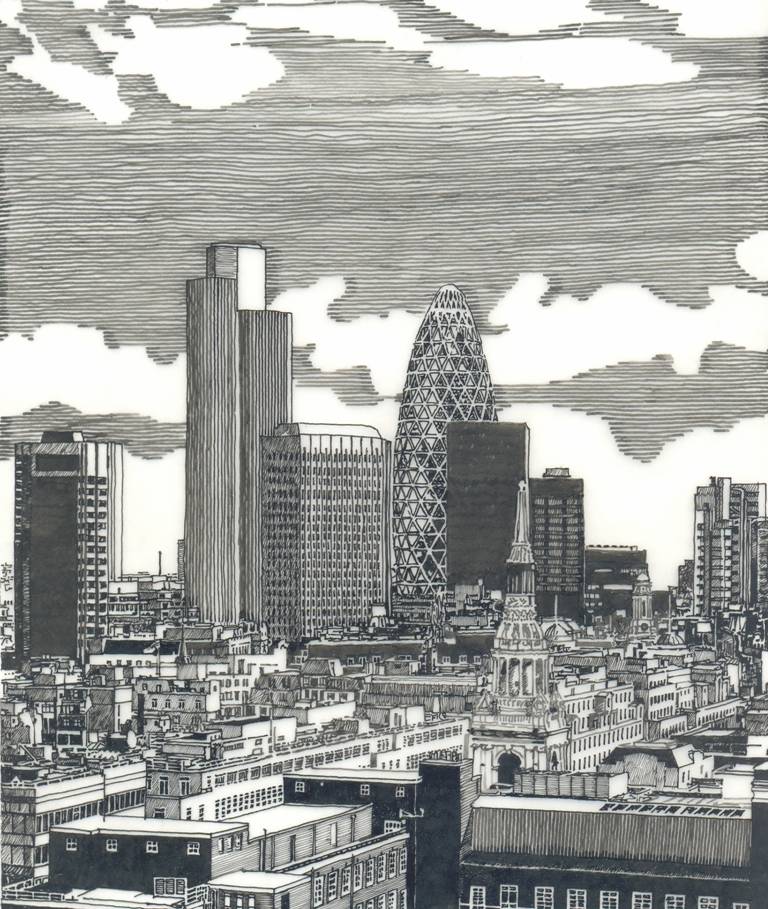}}
\caption{Architect: Norman Foster}
\label{fig5}
\end{center}
\end{figure}

When the time is not taken into account, we have what we call steady cells 
\cite{cc06,cc09}, defined in the 3D Cartesian orthogonal space O$xyz$.
 For instance, the half-sphere is an example of an Euclidean form
 (Figure \ref{fig5}:a), while an example of a non-Euclidean form 
(Figure \ref{fig5}:b) can be the set defined by
\begin{equation}\label{set}
\{(x,y,z)\in A\times \mathbb{R}:\ z=H-b(x^2+y^2)\},
\end{equation}
where
\[ 
A=\{(x,y)\in \mathbb{R}^2:\ x^2+y^2\leq H/b\},
\]
with $H$ denoting the height and 
$b$ is a constant depending on the radius of the circumference at the ground level. This obtained form is regular and aerodynamic such that has two principal features: it redistributes the gravity and deflects the gusting winds that often exist around the skyscrapers \cite{Ril}.

\begin{figure}
\begin{center}
\resizebox{1\hsize}{!}{\includegraphics{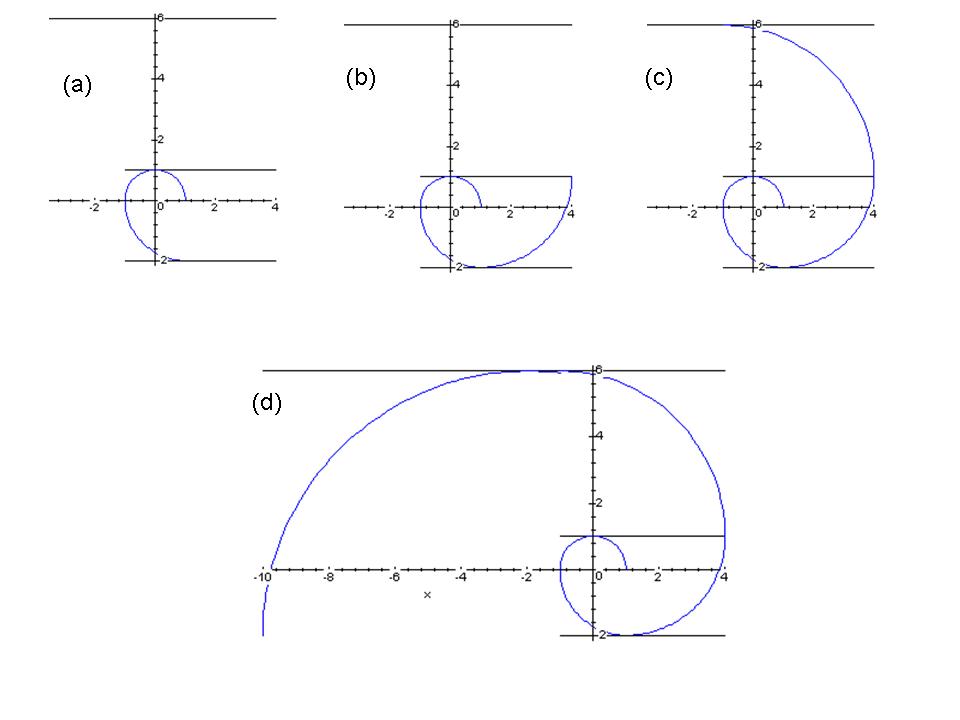}}
\caption{The construction of the Fibonacci spiral for squares 
with edge lengths (a) 1, 1, 2; (b) 3; (c) 5; (d) 8.}
\label{fig6}
\end{center}
\end{figure}

\begin{figure}
\begin{center}
\resizebox{1\hsize}{!}{\includegraphics{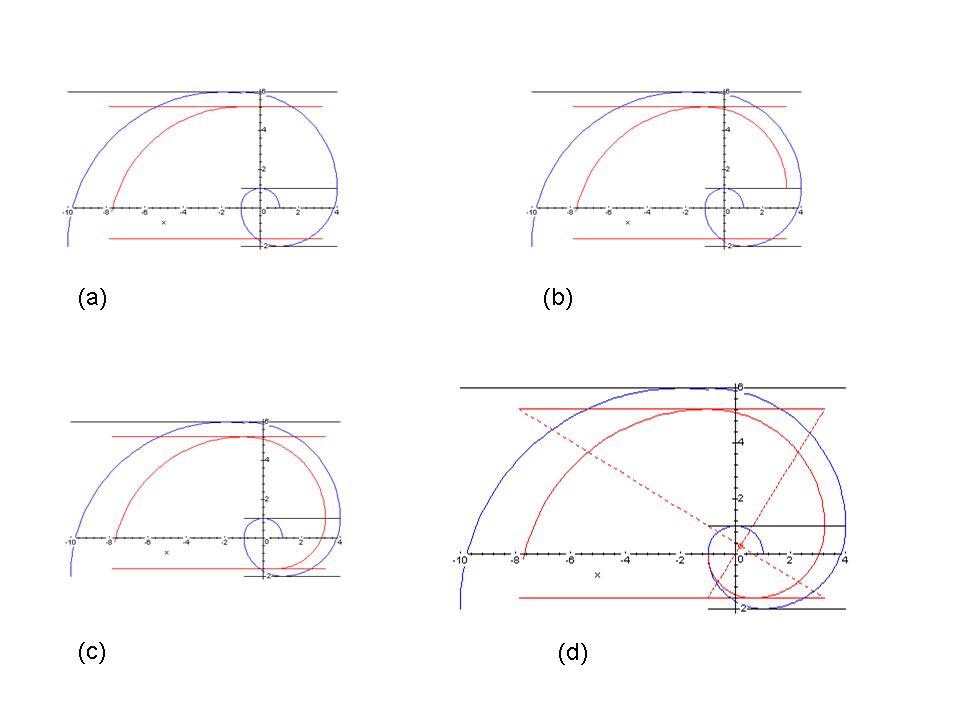}}
\caption{The Fibonacci spiral in blue (squares under 
the sequence 1, 1, 2, 3, 5, 8, $\cdots$), and  the construction of
 Divina spiral in red (under the golden rectangle).}
\label{fig7}
\end{center}
\end{figure}

\begin{figure}
\begin{center}
\subfigure[]{\includegraphics[width=5cm]{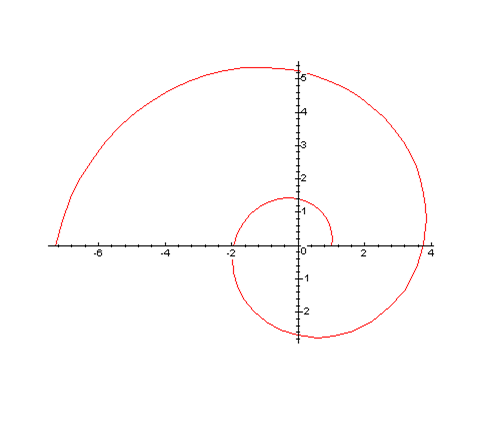}}
\subfigure[]{\includegraphics[width=5cm]{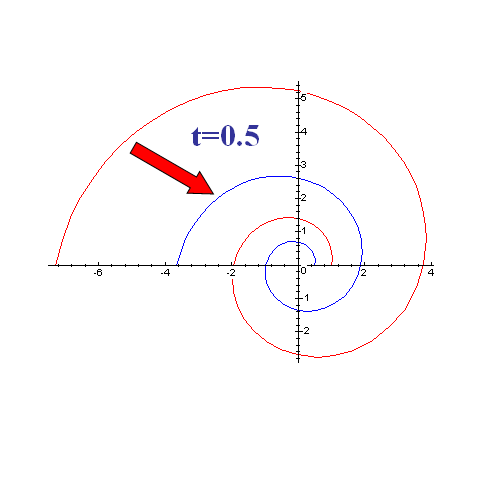}}
\caption{(a) The Spira mirabilis of Bernoulli, where the radius 
$r$ is defined by the relation $r=\phi^{b\theta}$ with $\theta$
 denoting the angle; (b) its space-time concept.}
\label{fig8}
\end{center}
\end{figure}

Next we discuss the real meaning of growth. It is known that a polygon is generated by overlapping suitable plane figures. The Fibonacci spiral is drawn with compass and straightedge 
(see Figure \ref{fig6}). At each iteration, the arc is similar to the next one. The outer part grows at a constant rate but faster than the inner part. The golden spiral is also a parametric variation of semi-circumferences 
(see Figure \ref{fig7}). Thus the classical form is a sequence of iterated motifs. The movement of a segment of a straight line around a fixed point originates a circle. Then the form is movement from a source. However these movements are spatial not attributing the time variable. Our objective is to illustrate a continuous spiral does not grow in its original meaning but also it can change in size as a function depending on space and time. 
A spiral can be defined by specifying the distance from the
origin or radius for each given angle, i.e., by a specific relation
between the radius and the angle.  
Considering the spiral 
(Figure \ref{fig8}:a), where its specific relation can be understood as
 the radius is a function on the angle or equivalently the radius satisfies an ordinary differential equation:
\[
r=\phi^{b\theta}\Leftrightarrow
{dr\over d\theta}=b\ln(\phi)r,
\]
 it can indeed be defined by the set, for every time parameter $t>0$,
\[
\{(x,y)\in \mathbb{R}^2:\ \phi^{b\arctan [y/x]}=\sqrt{x^2+y^2}
\},
\]
where $\phi$ denotes the golden number and 
$b$ is some constant. Then the continuous growth appears, if we take 
$t=1$ for the Spira mirabilis in red then it has its transformation for 
$t=0.5$ in blue 
(Figure \ref{fig8}:b). The idea that the form must have Euclidean characteristics is given up, i.e. the architectural object is obtained by the graphical representation of the mathematical mappings. The effects of selective growth capable of producing architecturally meaningful forms can be computationally simulated \cite{rou}.

\begin{figure}
\begin{center}
\subfigure[]{\includegraphics[width=5cm]{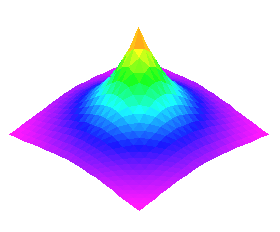}}
\subfigure[]{\includegraphics[width=5cm]{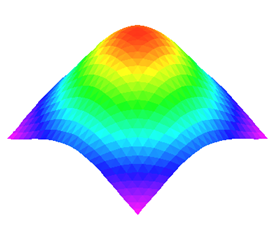}}
\subfigure[]{\includegraphics[width=5cm]{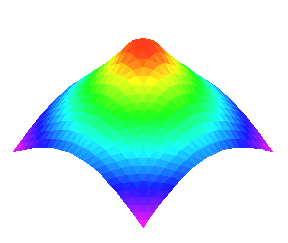}}
\caption{The graphical representation of a form defined by the set
$\{(x,y,z)\in [-2,2]^2\times\mathbb{R}:$ $z=\exp[-(x^2+y^2)^{1/t}]\}$,
 for the time parameter $t>0$: (a) $t=1$; (b) $t=2$; (c)
The graphical representation of a (steady-state) form defined by the set 
$\{(x,y,z)\in [-2,2]^2\times\mathbb{R}:$ $z=-(x^2+y^2)\sin
[1/\sqrt{x^2+y^2}]\}$.}
\label{fig12}
\end{center}
\end{figure}

Although the scientific community acknowledge the beauty of the mathematics as being universal, the aesthetic value of the form should be always given by the intensity of sensations (known as the emotion of the form 
\cite{cc09}), and never by the  mathematical beauty. 
Using mathematics without taken care of the expression, we obtain simply different shelters.
The distinction between the Euclidean and non-Euclidean objects, and their
correspondence to forms, disappears.
From the differential geometric point of view,
 the half-sphere (ex. of an Euclidean form) and the set defined by (\ref{set})
 (ex. of a non-Euclidean form) are equivalent (diffeomorphic),
 and even equivalent to a flat disc, say.
We emphasize that from the topologic point of view, also the forms illustrated
in Figure \ref{fig12} are equivalent (homeomorphic) to the above ones.
The distinction lies at the root of expression established in Section \ref{ef}. Notwithstanding, the expression of V-shaped cut given by the singular (non-differentiable) point at the top is the first inroad to rupture the way of thinking (Figure \ref{fig12}). 
 
\section{Conclusions}

The analysis of the form in the sphere of post-modernist art theory does not agree with the classic theory on the harmony of the proportions. Synthetically, the form into its components, expression and content, results from: 1) the properties of continuity and differentiability given at the 3D space, 2) the movement character given by the time parameter, 3) the knowledge of Architecture History. Without this last knowledge, there is the risk of falling into ambiguous signs and consequently into an emptiness of a society without aesthetics.

Therefore, the mathematical mappings should have a primordial role in the elaboration of the object as a sign where the man lives and it promotes an emotion, and contribute to one civilization more friendly to the environment. The research of living and aesthetic qualities demands a profounder knowledge of the new materials and of the capacity of their use.  
In conclusion, the relevance of our mathematical modeling  relies on the
 morphology that answers both to the human and earth needs. 
 Although the space-time, the growthness, and the mathematical language are
 beginning to be well established programmatic ideas since many years
 in the area of design in general \cite{leg},
 the grammar of composition gains a new concept: morphocontinuity,
 where  form is mathematically formulated and biophysically
 reinterpreted.

\section*{Acknowledgement}
The authors gratefully thank to Teot\'onio Agostinho for Figs. 1-6, and 8.

\end{document}